\documentclass[12pt]{article}
\usepackage{amssymb}
\makeatletter
\date{}

\newtheorem{theorem}{Theorem}[section]

\newtheorem{proposition}[theorem]{Proposition}

\newtheorem{corollary}[theorem]{Corollary}

\newtheorem{problem}[theorem]{Problem}

\newcommand{\edim}{{\rm e}$-${\dim}}

\newcommand{\z}{{\Bbb Z}}
\newcommand{\q}{{\Bbb Q}}

\newcommand{\N}{{\Bbb N}}

\newcommand{\lo}{\longrightarrow}

\newcommand{\black}{{\blacksquare}}

\begin{document}

\title{Resolving rational cohomological dimension
via \\ a Cantor group action}

\author{  Michael Levin\footnote{the author was supported by ISF grant 836/08}}

\maketitle
\begin{abstract} 
By  a Cantor group we mean a topological group homeomorphic 
to the Cantor set. We show that a compact metric space of
rational cohomological dimension $n$ can be obtained as 
the orbit space of a Cantor group action on a metric 
compact space of covering dimension $n$. Moreover,
the action  can be assumed to be free if $n=1$.
\\\\
{\bf Keywords:}  Cohomological Dimension,  Transformation Groups
\bigskip
\\
{\bf Math. Subj. Class.:}  55M10, 22C05 (54F45)
\end{abstract}
\begin{section}{Introduction}\label{intro}
Throughout this paper we assume that  maps are continuous and
spaces are separable metrizable. We recall that a  compactum
means a compact metric space. By the dimension $\dim X$  of a space $X$ we
assume the covering dimension.

Let  $G$ be  an abelian group. The cohomological dimension $\dim_G X $ of a space $X$
is the smallest integer $n$ such that  the Cech cohomology $H^{n+1} (X, A; G)$ 
vanishes for every closed subset $A$ of $X$.
Clearly $\dim_G X \leq \dim X$ for every abelian group $G$. Exploring connections
between cohomological and covering dimensions is one of the central topics
in Dimension Theory. Let us mention a few important results in this direction.
By the classical result of  Alexandroff 
$\dim X =\dim_\z X$ if $X$ is finite dimensional. Solving
a long standing open problem Dranishnikov constructed in
\cite{dranish-dim-3} an infinite dimensional
compactum $X$ with $\dim_\z X=3$. Edwards'
 famous cell-like resolution
theorem   \cite{edwards, walsh}
 asserts that every compactum of $\dim_\z=n$ is  the image
of an $n$-dimensional compactum under a cell-like map.  A map is cell-like 
if its fibers are cell-like compacta  and a compactum  is cell-like
if any map from it to a CW-complex is null-homotopic. 
 Edwards' 
theorem was extended in \cite{levin-rational} to rational
cohomological dimension: every compactum of $\dim_\q =n, n\geq 2,$ 
is the image of an $n$-dimensional compactum under a rationally acyclic map.
A map is rationally acyclic if its fibers are rationally acyclic compacta
and  a compactum  is rationally acyclic if its  reduced  Cech cohomology
with rational coefficients 
vanishes. Note that by the Begle-Vietoris theorem a cell-like map
and a rationally acyclic map cannot raise the integral and the rational cohomological
dimensions respectively.
This paper is devoted to establishing another connection
between rational cohomological and covering dimensions.

\begin{theorem}
\label{main-theorem}
Let $X$ be a compactum with $\dim_\q X =n$. Then there is an $n$-dimensional
 compactum  $Z$
and an action of a Cantor group  $\Gamma$ on $Z$ such that $X=Z/\Gamma$.
Moreover, the action of $\Gamma$ can be assumed to be free if $n=1$.
\end{theorem}

By a Cantor group we mean a topological group homeomorphic to the Cantor set.
Since a Cantor group is a pro-finite group one can easily derive from 
basic properties of transformation groups  \cite{bredon} 
that for an action of a Cantor group $\Gamma$ on a compactum $Z$ we have
$\dim_\q Z/\Gamma \leq \dim_\q Z$ and hence $\dim_\q Z/\Gamma \leq \dim Z$.
On the other hand  Dranishnikov and Uspenskij
 \cite{d-u} showed that:

\begin{theorem}
{\rm (\cite{d-u})}
\label{dranishnikov-uspenskij-theorem}	
	Let $f : Z  \lo X$ be a 0-dimensional map  of compacta $Z$ and $X$.
	Then $\dim_G Z \leq \dim_G X $ for every abelian group $G$. 
 
\end{theorem}	
Thus 
 for an action of a Cantor group
$\Gamma$ on a compactum $Z$ we have that
$\dim_\q Z/\Gamma =\dim_\q Z$ and hence Theorem \ref{main-theorem} provides a characterization of rational
cohomological dimension in terms of  Cantor groups actions.

Let us also note that, in general, the action of $\Gamma$ in Theorem \ref{main-theorem}
cannot be free for $n>1$. For example, consider 
a $4$-dimensional compact absolute retract $X$ with $\dim_\q X \leq 3$ 
constructed by Dranishnikov \cite{dranish-ANR, dranish-atlas}
and assume that $X$ is the orbit space of a free action of a Cantor
group $\Gamma$ on a compactum $Z$. One can easily observe that then
$Z=X \times \Gamma$ and hence $\dim Z = \dim X =4$.

A result closely related to Theorem \ref{main-theorem} can be derived
from the work of Dranishnikov and West \cite{dranish-west}. For a   prime
 $p$ we denote  by  $\z_p=\z/p\z$ the $p$-cyclic group
and by  $\z_p^\N=\Pi_{i=1}^\infty (\z_p)_i$ the product  of countably many 
copies of $\z_p$. It was  shown in \cite{dranish-west} that:
\begin{theorem}
{\rm (\cite{dranish-west})}
\label{theorem-dranishnikov-west}
Let $X$ be a compactum. Then for every prime $p$
there is a compactum $Y$ and an action of the group 
$\Gamma=\z_p^\N$ on $Y$ such that
$\dim_{\z_p} Y\leq 1$ and $X =Y/\Gamma$.
\end{theorem}

Consider any compactum $X$. By Theorem \ref{theorem-dranishnikov-west}
for every prime $p$ there is a compactum $Y_p$ and an action
of $\Gamma_p=\z^\N_p$ on $Y_p$ such that and $\dim_{\z_p} Y_p \leq 1$
and $X=Y_p/\Gamma_p$. Let ${\gamma}_p : Y_p \lo X$
be the projection, 
$\cal P$ the set of prime numbers and 
$Y=\Pi_{p \in \cal P} Y_p$.  Denote by $Z $   the pull-back 
of the maps $\gamma_p , p\in \cal P$, which is  the subset  of $Y$
consisting of the points $y \in Y$ such that $\gamma_p(y)=\gamma_q(y)$
for every $p,q \in \cal P$. Consider the group 
$\Gamma=\Pi_{p \in \cal P} \z^\N_p$ and the pull-back
action of $\Gamma$ on $Z$ defined by
$gy=(g_p y_p)$ for $g=(g_p)\in \Gamma$
and $y=(y_p) \in Z$  and notice that $X=Z/\Gamma$. 
Since the orbits
of $\Gamma_p$ on $Y_p$ are $0$-dimensional we get that the orbits of 
$\Gamma$ on $Z$ are $0$-dimensional and  the projection of $Z$
to $Y_p$ is $0$-dimensional for every $p$. 
 	By Theorem \ref{dranishnikov-uspenskij-theorem} we get
	that  $\dim_\q Z \leq \dim_\q X$ and  $\dim_{\z_p} Z \leq 1$
	for every $p \in \cal P$.  Then, 
	by Bockstein inequalities \cite{kuzminov-bockstein,dranish-atlas}, 
	$\dim_\z Z \leq \max \{\dim_\q  X, 2\}$  and we obtain
	
	\begin{theorem}{\rm (derived  from \cite{dranish-west})}
	\label{derived-dranishnikov-west}
	Let $X$ be a compactum with $\dim_\q X =n$. Then there is
	a compactum $Z$  and an action of the  group $\Gamma=\Pi_{p \in \cal P} \z^\N_p$
	on $Z$ such that $\dim_\z Z \leq \max\{n,2\}$ and $X =Z/\Gamma$.
	\end{theorem}

Theorem \ref{derived-dranishnikov-west} motivates
\begin{problem}
\label{problem}
Can the group $\Gamma$ in Theorem \ref{main-theorem} be assumed to be 
abelian? the product of countably many finite  cyclic groups?
\end{problem}

Let us finally mention that the interest in Cantor groups actions
is inspired by the Hilbert-Smith conjecture claiming  that
no Cantor group  effectively acts on a  manifold.
Smith \cite{smith} reduced
 the conjecture to the actions of 
 the groups of $p$-adic integers $A_p$, $p \in \cal P$, and
Yang \cite{yang} showed that if $A_p$ acts effectively on
a   manifold $M$ then $\dim_\z M/A_p =\dim M +2$ and if 
$A_p$ acts on a finite dimensional compactum $Z$ then
$\dim_\z Z/A_p \leq \dim Z +3$. Theorem \ref{main-theorem}
 shows that the dimensional restrictions imposed by $A_p$
do  not apply to  general  Cantor groups even in the following
extreme form: there is a free  action of a Cantor group on
a one-dimensional compactum raising  the integral 
cohomological dimension of
the orbit space to infinity. Indeed, take an infinite dimensional
compactum $X$ with $\dim_\q X=1$ and $\dim_\z X=\infty$.
 Then Theorem \ref{main-theorem}
produces a one dimensional compactum $Z$ and a free action of a Cantor
group $\Gamma$ on $Z$ such that $X=Z/\Gamma$. This example can be
considered as  complementary to 
 Dranishnikov-West's example \cite{dranish-west}
of  an action of $\z_p^\N$ on a $2$-dimensional compactum
raising the dimension of the orbit space to infinity.
\end{section}

\begin{section}{Preliminaries}

 Let us recall basic  definitions and results in Extension Theory and Cohomological
 Dimension that will be used in the proof
 of  Theorem \ref{main-theorem}.
 
 The extension dimension of a space $X$ is said to be  dominated by 
 a CW-complex $K$, written $\edim X \leq K$, if every map $f : A \lo K$ from
 a closed subset $A$ of $X$ extends over $X$.  Note
 that the property $\edim X \leq K$ depends only on the homotopy type
 of $K$.  
 The covering and cohomological dimensions can be 
 characterized by the following extension properties: 
 $\dim X \leq n$ if and only if the extension dimension of $X$ is  dominated by
 the $n$-dimensional sphere $S^n$ and 
 $\dim_G X \leq n$ if and only if the extension dimension of 
 $X$ is  dominated by the Eilenberg-Mac Lane complex
 $K(G,n)$. The extension dimension shares many properties of covering dimension.
 For example: 
 if $\edim X \leq K$ then for every $A \subset X$ we have
 $\edim A \leq K$, and if $X$ is a countable union of closed subsets whose
 extension dimension is dominated by $K$ then $\edim X \leq K$. 
  Let us list a few more properties.

 \begin{theorem}
 \label{dydak-union-theorem}
 {\rm \cite{dydak-union}}
 Let $X$ be a space and $K$ and $L$ CW-complexes.
 If   $X=A \cup B$ is the union of  subspaces $A$ and $B$ such that
 $\edim A \leq K$ and $\edim A \leq L$ then 
 the extension dimension of $X$ is  dominated
 by the join $K *L$.
 \end{theorem}
 \begin{theorem}
 \label{dranishnikov-splitting-theorem}
{\rm  \cite{dranishnikov-splitting}}
 Let $K$ and $L$ be countable CW-complexes and  $X$ a compactum such that
 $\edim X \leq K*L$. Then $X$ decomposes into subspaces $X=A \cup B$ 
 such that $\edim A  \leq K$ and $\edim B \leq L$.
 \end{theorem}

\begin{proposition}
{\rm (\cite{cell-like-levin})}
\label{extension-proposition}
Let $X$ be a compactum and $K$ a simply connected CW-complex such that
$K$ has only finitely many non-trivial homotopy groups and
$\dim_{\pi_i(K)} X \leq n$ for every $i>0$. Then $\edim X \leq K$.
\end{proposition}

\begin{theorem}
{\rm ( \cite{dranishnikov-extension})}
\label{dranishnikov-extension-theorem}
Let $X$ be a compactum and $K$ a CW-complex such that
$\edim X \leq K$. Then $\dim_{H_n(K)} X \leq n$ for every $n>0$.
\end{theorem} 
By a Moore space $M(\q,n)$ we will mean the model which is
the infinite
telescope of a sequence of maps from $S^n$ to $S^n$ of
all possible non-zero degrees. Note that $M(\q,1)=K(\q,1)$.

\begin{proposition}
\label{classifying}
Let $X$ be a compactum. Then $\dim_\q X \leq n$
if and only if $\edim X \leq  M(\q,n)$.
\end{proposition}
{\bf Proof.} By Theorem \ref{dranishnikov-extension-theorem}
the condition $\edim X \leq M(\q,n)$ implies $\dim_\q X \leq n$.
Let us show  that $\dim_\q X \leq n$  implies $\edim X \leq M(\q,n)$.
The case $n=1$ is obvious since $M(\q,1)=K(\q,1)$. Assume that $n>1$.
Then $M(\q, n)$ is simply connected and
 we have $H_i(M(\q,n))=0$ if $0< i <n$ and  
$H_i(M(\q, n))\otimes \q=H_i((M(\q,n))$ if $i \geq n$
  and hence
we have $\pi_i(M(\q,n))=0$ if $0< i <n$  and 
 $\pi_i(M(\q,n)) \otimes \q =\pi_i(M(\q,n))$ if $i \geq n$
 \cite{rational-book}.
Thus, by Bockstein Theorem \cite{kuzminov-bockstein, dranish-atlas}, $\dim_{\pi_i(M(\q,n))}  X \leq \dim_\q X \leq n$
for every $i \geq n$.  
Note that $M(\q,n)$ is the direct limit of its finite subtelescopes which
are homotopy equivalent to $S^n$. Then,
since 
$\pi_i(S^n)$ is   finite for $i\geq 2n$,  we have that $\pi_i(M(\q, n))$
is torsion for $i \geq 2n$ and hence $\pi_i(M(\q, n))=\pi_i(M(\q, n)) \otimes \q=0$
for $i \geq 2n$. Thus, by Proposition \ref{extension-proposition},
$\edim X \leq M(\q,n)$ and we are done. $\black$

\begin{corollary}
\label{corollary}
Let $X$ be a compactum. Then  $\dim_\q X \leq n$, $n>1$ if and only if
 $X$ decomposes into $X=A \cup B$ such that
$\dim_\q A \leq 1$ and $\dim B \leq n-2$.
\end{corollary}
{\bf Proof.} Note that $M(\q,n)$ is homotopy equivalent
to $\Sigma^{n-1} M(\q,1)=S^{n-2} * M(\q,1)$ and the corollary follows
from Proposition \ref{classifying} and Theorems \ref{dydak-union-theorem} and
\ref{dranishnikov-splitting-theorem}. $\black$

\begin{proposition}
\label{rational-dimension-via-spheres}
Let $X$ be 
a compactum. Then $\dim_\q X \leq n$ if and only if
for every map $f: A \lo S^n$ from a closed subset
$A$ of $X$ to a sphere $S^n$ we have that $f$ followed by 
a map of non-zero degree from $S^n$ to $S^n$ extends
over $X$.
\end{proposition}
{\bf Proof.} 
 The proof follows from Proposition \ref{classifying}
 and an easy adjustment of the proof of Proposition 3.13 of \cite{exotic-first}.
 $\black$
 \\\\
By a partial map from a compactum $X$ to a CW-complex  $K$ 
we mean a map $f : F \lo K$ from a closed subset $F$  of $X$ 
to $K$. A collection ${\cal F}$ of partial maps from $X$ to $K$
is said to be representative if 
for every partial map $f' : F' \lo K$  there is 
a map $f : F \lo K$ in ${\cal F}$ such that 
$F' \subset F$ and $f$ restricted to $F'$ is homotopic to $f'$.
Let $X$ be the inverse limit of a sequence of compacta $X_n$
with bonding maps $\omega_{n+1} : X_{n+1} \lo X_n$. 
We say that a partial map $f : F \lo K$ from $X_i$ to $K$ is extendable
in the inverse system if there is $j >i$ such that
for the map 
$\omega^i_j=\omega_{i+1} \circ \omega_{i+2} \circ \dots \circ \omega_j
: X_j \lo X_i$ we have that
$\omega^i_j$ restricted to $(\omega^i_j)^{-1}(F)$ and followed by
$f$ extends over $X_j$. The proof of the following
 proposition is simple and left to the reader.
\begin{proposition}
\label{representative-extendable}
Let $X$ be a compactum and $K$ a CW-complex.

(i) If $K$ is a countable CW-complex  then there is
a countable representative collection of partial maps from $X$ 
to $K$;

(ii) Let $X$ be the inverse limit of compacta $X_n$ and let
${\cal F}_n$ be a representative collection
of partial maps from $X_n$ to $K$ such that
for every $n$ and every $f$ in ${\cal F}_n$ we have
that $f$ is extendable in the inverse system. Then $\edim X \leq K$.
\end{proposition}

\end{section}
\begin{section}{Proof of Theorem \ref{main-theorem}}

\begin{proposition}
\label{extension-via-covering}
Let $(K, L)$ be a  pair of finite simplicial complexes with $K$ being
connected, 
$\psi : S^1 \lo S^1$ a  covering map
 and $\phi : L \lo S^1$ a map such that
$\phi$ followed by $\psi$ extends to a map $g : K \lo S^1$.
Consider  a component $M$ of   the pull-back  of the maps $g$ and $\psi$
and the projection $\beta_M : M \lo K$.  Then 
$\beta_M$ restricted to $\beta_M^{-1}(L)$ and followed by $\phi$
extends over $M$.
\end{proposition}
{\bf Proof.}  
Clearly $\beta_M$ is a covering map. 
Let $m=\deg \psi$ and assume that  $\z_m =\z/m\z$
 freely  acts on $S^1$ so that 
$\psi$ is  the projection  to the orbit space $S^1=S^1/\z_m$.
Then 
the action of $\z_m$ on $S^1$ induces the corresponding
free action of $\z_m$ on the pull-back $Z$  of the maps $g$ and $\psi$.

Consider a component $L' $ of $L$ and let $N'$ be a component of
$\beta_M^{-1}(L')$. 
Clearly  $\beta_M$ restricted to $N'$ and $L'$ is
a  covering map. 
Let us show that $\beta_M$ is one-to-one on $N'$.
Take
a map $\alpha : S^1 \lo L'$. Note that  $\phi \circ \alpha$ is a  lifting 
of  $\psi\circ \phi  \circ \alpha$ 
via $\psi$ and hence the map 
$(\alpha, \phi \circ \alpha) : S^1 \lo Z\subset K \times S^1$
is a lifting of $\alpha$
via the projection of $Z$ to $K$. Then  $(\alpha, \phi \circ \alpha)$
followed
by (the action of) an element of $\z_m$ provides a lifting 
of
$\alpha$  to $N'$ and hence $\beta_M$ is one-to-one
on $N'$. 

Let $\rho : M \lo S^1$ be the  projection  to $S^1$.
Note that $\rho$ restricted to $N'$ is a lifting
 (via the map $\psi$)
 of the $\beta_M$ restricted
to $N'$ and followed by $\psi \circ \phi$.
Then the maps $\rho|N'$ 
and $\phi \circ (\beta_M |N')$ coincide up to 
the action of $\z_m$  on $S^1$ and hence
they are homotopic. Thus $\rho$ extends up to homotopy
the map $\beta_M$ restricted to $N'$ and followed by $\phi$  and
hence $\beta_M$ restricted to $\beta_M^{-1}(L)$ and followed by
$\phi$ extends over $M$. The proposition is proved. $\black$
\\

Let finite groups $\Gamma$, $\Gamma'$ act on
compacta $Y$ and $Y'$ respectively and let 
$\alpha : \Gamma' \lo \Gamma $ be a homomorphism
and $\beta : Y' \lo Y$ a map.  We say that $\alpha$ and $\beta$
agree  if 
$\beta(gy)=\alpha(g)\beta(y)$
for every $g \in \Gamma'$ and $y \in Y'$. 

 \begin{proposition}
 \label{auxiliary-proposition}
 Let $X$ and $Y$ be   compacta, $A$ a closed connected  subset
 of $X$ with $\dim_\q A \leq 1$,
 $\Gamma$   a finite group acting on
 a compactum $Y$ with  $X=Y/\Gamma$ such that
 $\gamma^{-1}(A)$
 is connected  and 
 $\Gamma$ freely acts 
  on $ \gamma^{-1}(A)$
where 
$\gamma : Y \lo X=Y/\Gamma$
is the projection.  Then
for every map $f: F \lo S^1$ from a closed subset
$F$ of 
$\gamma^{-1}(A)$ to a circle $S^1$ there are a compactum $Y'$,
 a finite group  $\Gamma'$ acting on $Y'$, an epimorphism 
 $\alpha : \Gamma' \lo \Gamma$ and
  a surjective  map $\beta : Y' \lo Y$
 so that $\alpha$ and $\beta$ agree,
  $X=Y'/\Gamma'$ 
and for the projection
 $\gamma' : Y' \lo X=Y'/\Gamma'$ 
  we have that $\gamma' = \gamma \circ \beta$,
  $\gamma'^{-1}( A)$ is connected,
 $\Gamma'$ freely acts on $\gamma'^{-1}( A)$ and $\beta$ restricted
 to $(\beta)^{-1}(F)$ and followed by $f$
  extends over $\gamma'^{-1}(A)$.

 \end{proposition}
 {\bf Proof.} Clearly $\dim_\q \gamma^{-1}(A) \leq 1$ and hence,
 by Proposition \ref{rational-dimension-via-spheres},
 there is a covering map $\psi : S^1 \lo S^1$ so that
 $\psi \circ f$ extends to over $ \gamma^{-1}(A)$.
 Since $\Gamma$ is finite,
 the action of $\Gamma$ is free  on $\gamma^{-1}(A)$ and
  $\gamma^{-1}(A)$ is compact    and connected one can approximate 
 the action
 of $\Gamma$ over an closed  neighborhood  $B$ of $A$ 
 through a free action of $\Gamma$ on a finite connected simplicial 
 complex $K$  and a     surjective  map 
  $\mu_Z : Z=\gamma^{-1}(B) \lo K$ such that $\mu_Z$ 
	commutes with the actions of 
	$\Gamma$  and $f$ factors
  up to homotopy through $\mu_Z$ restricted to $F$ and
	a map $\phi : L \lo S^1$ from
	a subcomplex $L$ of $K$ so that $\mu_Z(F) \subset L$  and
	$\phi$ followed by $\psi$ extends to a map  $g : K \lo S^1$.
	Let 
	$\mu : B=Z/\Gamma \lo K/\Gamma$ be the map induced by $\mu_Z$
and 
 $\gamma_K  : K \lo K/\Gamma$ the covering projection.

  In this proof we always consider a covering map as a pointed map
  and identify the fundamental group of
  the covering space with the subgroup
  of the fundamental group of the base space obtained
  under the induced monomorphism of the fundamental groups.
  Thus
  $\pi_1(K)$ is a normal subgroup  of $\pi_1(K/\Gamma)$
  with $\pi_1(K/\Gamma)/\pi_1(K)=\Gamma$. Let 
  $M$ be a connected component of
  the pull-back of $g$ and $\psi$ and let $\beta_M : M \lo K$ be the projection.
	Clearly $\beta_M$ is 
  a covering map.  Then
   $\pi_1(M)$ is a subgroup of $\pi_1(K)$
  of finite index and hence there is  a normal subgroup
  $G$ of $\pi_1(K/\Gamma)$ of finite index such that
  $G$ is contained in $\pi_1(M)$.
  Consider   covering maps $\gamma'_K : K' \lo K/\Gamma$ and 
  $\beta_K : K' \lo K$
  from a covering space $K'$ with $\pi_1(K')=G$.
  Then $\beta_K$ factors through $\beta_M$ and 
  $\Gamma'=\pi_1(K/\Gamma)/G$ acts on $K'$ so that
  $K'/\Gamma'=K/\Gamma$ and 
  for the induced epimorphism $\alpha : \Gamma' \lo \Gamma$
  we have that $\beta_K$ and $\alpha$ agree.
  Denote by $Z'$ the pull back of $\gamma'_K$ and $\mu$,
  consider the pull-back action of $\Gamma'$ on $Z'$
  and let $\gamma'_Z : Z' \lo B= Z/\Gamma=Z'/\Gamma'$
  and $\beta_Z : Z' \lo Z$  be the projections induced 
  by $\gamma'_K$ and $\beta_K$ respectively.
	
	Take an open neighborhood $V $ of $A$ in $X$ so that $V \subset B$
	and denote $U=\gamma^{-1}(V)$ and $U'=\gamma'^{-1}_Z(V)=\beta_Z^{-1}(U)$.
	Set $Y'$ to  be the disjoint union of
	$Y\setminus  U$ and $U'$  and $\beta : Y' \lo Y$ to be
	the function 
	defined by the identity map on  $Y\setminus U$
	and by the map $\beta_Z$ on $U'$. Turn $Y'$ into a compactum
	by 
	preserving the topologies of $Y\setminus U$ and
	$U'$  and declaring the function $\beta$ to be the  quotient  map. 
	Extend the action of  $\Gamma'$ on $U'$  over $Y'$ by setting
	 $gy = \alpha(g) y$ for $y \in Y' \setminus U'=Y\setminus U$
	and $g \in \Gamma'$. It is easy to see that
		the action of $\Gamma'$ is well-defined,  
		$\alpha$ and $\beta$ agree,  $X=Y'/\Gamma'$ and
		for the projection $\gamma' : Y' \lo X=Y'/\Gamma'$ we have 
		$\gamma'=\gamma \circ \beta$. 
		
		By Proposition \ref{extension-via-covering} the map $\beta_M$
		restricted to $\beta_M^{-1}(L)$ and followed by $\phi$ extends
		over $M$. 
		Recall that $\beta_K$ factors through $\beta_M$ and hence
		$\beta_K$
		restricted to $\beta_K^{-1}(L)$ and followed by $\phi$ extends
		over $K'$. 
		Then $\beta_Z $ restricted to $\beta_Z^{-1}(F)$ and followed by
		$f$ extends over $Z'$  and hence
		 $\beta$ restricted to $\beta^{-1}(F)$ and followed by
		$f$ extends over $\beta^{-1}(A)$ and the proposition follows.
		$\black$

\begin{proposition}
\label{main-proposition}
Let $X$ be a compactum and $A$ a closed subset of $X$ with
$\dim_\q A \leq 1$. Then there is a compactum $Y$ and an action
of a Cantor group $\Gamma$ on $Y$ such that $X=Y/\Gamma$ and
for the projection $\gamma : Y \lo X$ we have that
$\dim \gamma^{-1}(A) \leq 1$ and the action of $\Gamma$ is free on $\gamma^{-1}(A)$.
\end{proposition}
{\bf Proof.} Note that without loss of generality
one can replace $X$ by any larger compactum and $A$ by any
larger closed subset  of $\dim_\q \leq 1$. Also note that
by adding to $A$ a set of $\dim \leq 1$ we still preserve 
$\dim_\q \leq 1$. Take a Cantor set $C$ 
embedded into an interval $I$
 and 
a surjective map $\phi : C \lo A$ and consider
the compactum obtained by attaching to $I$ the mapping 
cylinder of $\phi$. Then attaching this compactum to $X$
through the set $A$ we enlarge $A$ to a connected set by adding
a subset of $\dim =1$ and
hence we may assume that $A$ is connected.

Set $Y_1=X$ and $\Gamma_1 =\{ 1\}$.
Apply Proposition \ref{auxiliary-proposition} to
construct by induction  for every $n$ 
 a compactum  $Y_n$, a finite group  $\Gamma_n$ acting on $Y_n$,
 an epimorphism $\alpha_{n+1} : \Gamma_{n+1} \lo \Gamma_n$
 and a surjective map $\beta_{n+1} : Y_{n+1} \lo Y_n$ 
such that $\alpha_{n+1}$ and $\beta_{n+1}$ agree, 
 $X=Y_n/\Gamma_n$
and for the projections $\gamma_n : Y_n \lo X=Y_n/\Gamma_n$
 we have
that $\gamma_{n+1}=\gamma_n \circ \beta_{n+1}$,
$\Gamma_n$ freely acts on $\gamma^{-1}_n(A)$ and for
a map, that will be specified latter, $f_n: F_n \lo S^1$ from
a closed subset $F_n$ of $\gamma^{-1}_n(A)$ we have
that the map $\beta_{n+1}$ restricted
to $(\beta_{n+1})^{-1}(F_n)$  and followed by $f_n$ extends
over $\gamma^{-1}_{n+1}(A)$.

Denote $Y={\rm invlim}(Y_n, \beta_n)$ and
$\Gamma= {\rm invlim }(\Gamma_n, \alpha_n)$. 
Clearly $\Gamma$ is a compact $0$-dimensional group,
$\Gamma$ acts on $Y$ so that $X=Y/\Gamma$ and $\Gamma$ acts 
freely on $\gamma^{-1}(A)$ where $\gamma : Y \lo X=Y/\Gamma$
is the projection. Let us show that 
on each step of the construction the map  $f_n$ can be chosen
in a way that guarantees that $\dim \gamma^{-1}(A) \leq  1$.

Once $Y_i$ is constructed, choose, 
by (i) of  Proposition \ref{representative-extendable},
 a countable representative collection
${\cal F}_i$ of partial maps from 
$A_i=\gamma_i^{-1}(A)$ to $S^1$ and fix a surjection
 $\tau_i : \N \lo {\cal F}_i$.
 Take any bijection $\tau : \N \lo \N \times \N$
so that for every $n$ and  $\tau(n)=(i,j)$ we have $i \leq n$. Consider
the inductive step of the construction from $n$ to $n+1$
and  take  $f=\tau_i(j) \in {\cal F}_i$ with $(i,j)=\tau(n)$.
 Recall that $i \leq n$
and hence the collection ${\cal F}_i$ is already defined.
Denote 
$\beta^i_n = \beta_{n} \circ \dots \circ \beta_i : Y_n \lo Y_i$
and  $F_n=(\beta^i_n)^{-1}(F)$ 
and let  $f_n =f\circ \beta_n^i | F_n  : F_n \lo S^1$ be
the map we use on the inductive step of the construction.
Then $f_n$ is extendable in the inverse system 
$\gamma^{-1}(A)={\rm invlim}(A_n, \beta_n|A_n)$
and 
 hence $f$ is extendable in this inverse system 
  as well. Thus, by (ii) of
 Proposition \ref{representative-extendable}, 
 $\edim  \gamma^{-1}(A) \leq S^1$ and hence $\dim \gamma^{-1}(A) \leq 1$.

If $\Gamma$ is a finite group then replacing 
both
$Y$ and $\Gamma$ by  the products $Y\times C$ and $\Gamma \times C$
with any Cantor group $C$ we may assume that $\Gamma$ is a Cantor group
and 
the proposition follows.
$\black$
\\\\
{\bf Proof of Theorem \ref{main-theorem}.} The case $n=0$ is trivial.
The case $n =1$ follows from Proposition \ref{main-proposition}.
Assume that $n \geq 2$.  By Corollary \ref{corollary}
decompose $X$ into $X=A \cup B$ so that $\dim_\q A \leq 1$ and
$\dim B \leq n-2$. Enlarging $B$ to a $G_\delta$-subset of $\dim \leq n-2$
and replacing $A$ by the smaller set
$X \setminus B$ we may assume that $A$ is $\sigma$-compact.
Represent $A$ as the union $A=\cup_{i=1}^\infty A_i$ of
compacta $A_i$. By Proposition \ref{main-proposition} there is
a compactum $Y_i$ and a Cantor group $\Gamma_i$ acting on $Y_i$
such that $X=Y_i/\Gamma_i$ and for the projection 
$\gamma_i : Y_i \lo X$ we have that $\dim \gamma^{-1}_i(A_i) \leq 1$.
Let $Y=\Pi_{i=1}^\infty Y_i$ and let $Z\subset Y$ be the pull-back
of the maps $\gamma_i$. Then 
$\Gamma=\Pi_{i=1}^\infty \Gamma_i$  is a Cantor group and for 
the pull-back action  of $\Gamma$  on $Z$ we have
 that $X=Z/\Gamma$.  Let $\gamma : Z \lo X$ be the projection.
Note that $\gamma$ is a $0$-dimensional map.
Also note that for every $i$  
the map $\gamma$ factors through the $0$-dimensional projection
$p_i : Z \lo Y_i$ and the  map $\gamma_i : Y_i \lo  X$. Thus we have, by Hurewicz Theorem,  
that
$\dim \gamma^{-1}(B) \leq \dim B \leq n-2$ and, since  
$\gamma^{-1}(A_i)=p_i^{-1}(\gamma_i^{-1}(A_i))$, we also have
 $\dim \gamma^{-1}(A_i)  \leq \dim \gamma_i^{-1}(A_i) \leq 1$ and
hence $\dim \gamma^{-1}(A) \leq 1$. 
Thus $\dim Z \leq \dim \gamma^{-1}(A) +\dim \gamma^{-1}(B) +1=n$. Recall that
$\dim_\q Z= \dim_\q X=n$ and hence $\dim Z=n$. The theorem is proved.
$\black$

 \end{section}

Michael Levin\\
Department of Mathematics\\
Ben Gurion University of the Negev\\
P.O.B. 653\\
Be'er Sheva 84105, ISRAEL  \\
 mlevine@math.bgu.ac.il\\\\

\begin{thebibliography}{99}

\bibitem{bredon}
Bredon, Glen E.
Introduction to compact transformation groups. 
Pure and Applied Mathematics, 46. New York-London: Academic Press. 
XIII,459 p.  (1972).

\bibitem{dranish-dim-3}
Dranishnikov, A. N.
On a problem of P. S. Aleksandrov, Mat. Sb. (N.S.) 135(177) (1988), 
no. 4, 551–557, 560. 

\bibitem{dranish-ANR}
 Dranishnikov, A. N. Homological dimension theory. 
Russian Math. Surveys 43 (4) (1988), 11-63.


\bibitem{dranishnikov-extension}
Dranishnikov, A. N. Extension  of mappings into CW-complexes. (Russian)  
Mat. Sb.  182  (1991),  no. 9, 1300--1310;  translation in  Math. USSR-Sb.  
74  (1993),  
no. 1, 47-56.

\bibitem{dranishnikov-splitting}
Dranishnikov, A. N. On the mapping intersection problem.  
Pacific J. Math.  173  (1996),  no. 2, 403-412.

\bibitem{dranish-atlas}
Dranishnikov, A. N. Cohomological dimension theory of compact metric spaces, Topology Atlas invited contribution, http://at.yorku.ca/topology.taic.html, 
arXiv:math/0501523.

\bibitem{exotic-first}
 Dranishnikov A., Levin M. 
Dimension of the product and classical formulae of dimension theory, Trans. Amer. Math. Soc., to appear. arXiv:1210.2775


\bibitem{d-u}
 Dranishnikov, A.N.; Uspenskij, V.V. Light maps and extensional dimension, 
Topology Appl. 80
(1997) 91-99.



\bibitem{dranish-west} 
 A.N. Dranishnikov and J.E. West, Compact group actions that raise dimension to infinity, Topology and its Applications 80 (1997), 101-114.
	arXiv:math/0212329


\bibitem{dydak-union}
Dydak, Jerzy Cohomological dimension and metrizable spaces. II.  
Trans. Amer. Math. Soc.  348  (1996),  no. 4, 1647-1661.

\bibitem{edwards}
R. D. Edwards, A theorem and a question related to cohomological dimension and 
cell-like maps, Notices of the AMS, 25(1978), A-259.


\bibitem{kuzminov-bockstein}
 Kuzminov, V.I. Homological dimension theory, Russian Math. Surveys 23, 
issue 5 (1968), 1--45.

\bibitem{cell-like-levin}
Michael Levin,
Cell-like resolutions preserving cohomological dimensions. 
 Algebr. Geom. Topol. 3 (2003) 1277-1289. 

\bibitem{levin-rational}
Michael Levin,
Rational acyclic resolutions.  Algebr. Geom. Topol. 5 (2005) 219-235.

\bibitem{smith}
P.A. Smith, Transformations of finite period, III, 
Newman’s Theorem, Ann. of Math. 42 (1941), 446-458. 

\bibitem{walsh}
J. Walsh, Dimension, cohomological dimension, and cell-like mappings. Shape theory and geometric topology (Dubrovnik, 1981), pp. 105–118, Lecture Notes in Math., 870, Springer, Berlin-New York, 1981.

\bibitem{rational-book}
Félix, Yves; Halperin, Stephen; Thomas, Jean-Claude (2001), 
Rational homotopy theory, Graduate Texts in Mathematics 205, 
New York: Springer-Verlag, pp. xxxiv+535.
\bibitem{yang}
 C.T. Yang, p-Adic transformation groups, Mich. Math. J. 7 (1960), 201-218. 


\end{thebibliography}
\end{document}